\documentclass[11pt]{article}
\usepackage{amsmath}
\usepackage{color}
\usepackage{graphicx}
\usepackage{amsmath,amsfonts,amssymb,amsthm}
\usepackage{graphicx}
\usepackage{epsfig}
\usepackage{epstopdf}
\usepackage{color}
\usepackage{fullpage}
\usepackage{setspace}
\usepackage{enumerate}
\def\tto{\;{\lower 1pt \hbox{$\rightarrow$}}\kern -10pt
\hbox{\raise 2pt \hbox{$\rightarrow$}}\;}

\begin{document}
\begin{center}
\vspace*{0.3in} {\bf ON A CONVEX SET WITH NONDIFFERENTIABLE \\METRIC PROJECTION}\\[1ex]

Shyan S. Akmal\footnote{Fariborz Maseeh Department of Mathematics and Statistics, Portland State University, Portland, OR 97202, United States (Email: shyan.akmal@gmail.com).}, Nguyen Mau Nam\footnote{Fariborz Maseeh Department of Mathematics and Statistics, Portland State University, Portland, OR 97202, United States (Email: mau.nam.nguyen@pdx.edu).}, and J. J. P. Veerman \footnote{Fariborz Maseeh Department of Mathematics and Statistics, Portland State University, Portland, OR 97202, United States, and CCQCN, Dept of Physics, University of Crete, 71003 Heraklion, Greece (Email: veerman@pdx.edu).}
\end{center}
\small{{\bf Abstract.} A remarkable example of a nonempty closed convex set in the Euclidean plane for which the directional derivative of the metric projection mapping fails to exist was constructed by A. Shapiro. In this paper, we revisit and modify that construction to obtain a convex set with \emph{smooth boundary} which possesses the same property.\\[1ex]
\bf AMS subject classifications.} 49J53, 49J52, 90C31

\newtheorem{Theorem}{Theorem}[section]
\newtheorem{Proposition}[Theorem]{Proposition}
\newtheorem{Remark}[Theorem]{Remark}
\newtheorem{Lemma}[Theorem]{Lemma}
\newtheorem{Corollary}[Theorem]{Corollary}
\newtheorem{Definition}[Theorem]{Definition}
\newtheorem{Example}[Theorem]{Example}
\renewcommand{\theequation}{\thesection.\arabic{equation}}
\normalsize

\section{A Convex Set with Smooth Boundary}
\setcounter{equation}{0}

Define a strictly decreasing sequence of real numbers $\{\alpha_n\}_{n\in\mathbb{N}}\subset (0,\pi/2]$ with
\begin{equation}\label{c1}
\lim_{n\rightarrow \infty}\alpha_n = 0 \; \mbox{\rm and }\alpha_{n+1}\leq \dfrac{\alpha_{n}+\alpha_{n+2}}{2}\; \mbox{\rm for all }n\in \mathbb N.
\end{equation}
Now we identify $\mathbb{R}^2$ equipped with the Euclidean norm $\|\cdot\|$ with $\mathbb{C}$ and let $A_n=e^{i\alpha_n}$. A beautiful and surprisingly simple example of a nonempty closed convex set for which the directional derivative of the metric projection mapping fails to exist was constructed by A. Shapiro in \cite{Shap}. This set  is essentially the
convex hull $J$ of the collection of points $0$, $1$, and $\{A_n\}_{n\in\mathbb{N}}$. Note that this set does not have smooth boundary. More positive and negative results on the existence of directional derivatives to the metric projection mapping as well as applications to optimization can be found in \cite{Abat,FP,KruskalC,bor,OD,Shap94,Shap,OldShap} and the references therein.

\begin{figure}[!ht]
\begin{center}
\includegraphics[width=7cm]{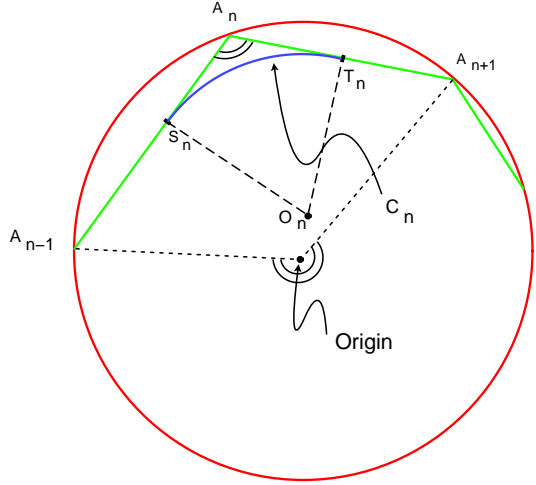}
\caption{The construction of a convex set with smooth boundary.}
\end{center}
\label{fig:convexset1}
\end{figure}

To define a convex set with smooth boundary, we start by choosing $\alpha_1=\pi/2$ and proceeding
as before to obtain the set $J$. The strategy to obtain a convex set $K$ with smooth boundary is to replace the pointy
parts of this figure by circular
arcs; see Figure \ref{fig:convexset1}. Let $T_n$ be the midpoint of the line
segment $A_nA_{n+1}$ and let $S_n$ the point in the line segment $A_{n-1}A_n$ so that
\begin{equation}
\|A_n-S_{n}\|= \|A_n-T_{n}\|= \sin\left(\dfrac{\alpha_n-\alpha_{n+1}}{2}\right).
\label{eq:An-Sn}
\end{equation}
Replace the two line segments $T_nA_n$ and $A_nS_n$ by a circular arc $C_n$ tangent to both
segments. Let $O_n$ be the center of the circle that contains $C_n$ as an arc and let $r_n$ denote
the radius of the circle. Let $J_1$ be the convex hull of the points $0$, $1$, the circular
arcs $\{C_n\}_{n\in\mathbb{N}}$, and the line segments connecting them. Let $J_2$ be the image of $J_1$ under reflection in the real axis and let
$J_3$ be the reflection of $J_1\cup J_2$ in the imaginary axis. Then we define $K:= J_1\cup J_2\cup J_3$. The set obtained has \emph{smooth boundary} in the sense we will define shortly.

\begin{Lemma}\label{lm1} $\lim_{n\rightarrow \infty}
\left| r_n -  2\,\dfrac{\alpha_n-\alpha_{n+1}}{\alpha_{n-1}-\alpha_{n+1}}\right| = 0$.
\end{Lemma}
\noindent{\bf Proof:} Consider the angle $\psi_n$ at $A_n$ and the angle $\phi_n$ and the origin as indicated
by the double arcs in Figure \ref{fig:convexset1}. From our definition of $\alpha_n$, we see that
$\phi_n=2\pi-(\alpha_{n-1}-\alpha_{n+1})$. By the Inscribed Angle Theorem, we have $\psi_n=\frac 12 \phi_n$. Thus,
\begin{equation}
\psi_n=\pi-\frac 12(\alpha_{n-1}-\alpha_{n+1}).
\label{eq:psi-n}
\end{equation}
The figure $A_nS_nO_nT_n$ is a \emph{right kite} with right angles at $S_n$ and at $T_n$. Therefore,
\begin{equation}
\dfrac{\|O_n-T_n\|}{\|A_n-T_n\|}=\tan \left(\frac{\psi_n}{2}\right) = \left[\tan{\left(\dfrac{\pi-\psi_n}{2}\right)}\right]^{-1}.
\label{eq:tangent}
\end{equation}
Using \eqref{eq:An-Sn}, \eqref{eq:psi-n}, and \eqref{eq:tangent}
in the relation
\begin{equation*}
r_n=\dfrac{\|O_n-T_n\|}{\|A_n-T_n\|}\,\|A_n-T_n\|,
\end{equation*}
we see that
\begin{equation}\label{formularn}
r_n=
\dfrac{\sin(\frac 12(\alpha_n-\alpha_{n+1}))}{\tan(\frac 14(\alpha_{n-1}-\alpha_{n+1}))}.
\end{equation}
The result then follows easily.
\qed

In what follows, we will distinguish three cases:\\
\hspace*{0.6in} Case $A$: $\alpha_n = C n^{-q}$, where $C, q>0$.\\
\hspace*{0.6in} Case $B$: $\alpha_n= C\lambda^n$, where $C>0$ and $\lambda \in (0,1)$. \\
\hspace*{0.6in} Case $C$: $\alpha_n= C\lambda^{n^2}$, where $C>0$ and $\lambda \in (0,1/2)$.

Recall that a subset $\Omega$ of $\Bbb R^m$ is called convex if
$$\alpha x+(1-\alpha)y\in \Omega\; \mbox{\rm whenever }x,y\in \Omega\; \mbox{\rm and }\alpha\in (0,1).$$
A function $f: \Omega\to \Bbb R$ is called convex if
$$f(\alpha x+(1-\alpha)y))\leq \alpha f(x)+(1-\alpha)f(y)\; \mbox{\rm for all }x,y\in \Omega\; \mbox{\rm and }\alpha\in (0,1).$$
If $-f$ is convex, then $f$ is called concave.

The lemma below will be important in what follows.

\begin{Lemma}\label{lemma12} In each of the three cases, the sequence $\{\alpha_n\}$ is strictly decreasing and satisfies the conditions in \eqref{c1}. Moreover, $\lim_{n\rightarrow \infty} r_n$ exists and:
$$
\lim_{n\rightarrow \infty} r_n= \left\{\begin{array}{cc} 1& \textrm{Case $A$}\\
                                                     \dfrac{2\lambda}{1+\lambda}& \textrm{Case $B$}\\
                                                     0& \textrm{Case $C$}
                                                     \end{array}\right.
$$
\end{Lemma}
\noindent{\bf Proof:} In Case A, it is obvious that $\{\alpha_n\}$ is strictly decreasing. Since the function $g(x):=x^{-q}$ is convex on $(0, \infty)$,
$$g(n+1)\leq \dfrac{g(n)+g(n+2)}{2},$$
which implies that $\alpha_{n+1}\leq \dfrac{\alpha_n+\alpha_{n+2}}{2}$.

By Lemma \ref{lm1}, we also see that
$$\lim_{n\to\infty}r_n=\lim_{n\to\infty}2\,\dfrac{\alpha_n-\alpha_{n+1}}{\alpha_{n-1}-\alpha_{n+1}}=2\lim_{n\to\infty}\dfrac{n^{-q}-(n+1)^{-q}}{(n-1)^{-q}-(n+1)^{-q}}=1.$$
The proof for Case B and Case C is left for the reader. \qed

\begin{Remark}{\rm In Case C, we can replace $\lambda\in (0,1/2)$ by $\lambda\in (0,1)$ and show that $\{\alpha_n\}$ satisfies conditions in \eqref{c1} for all $n\geq n(\lambda)$, where $n(\lambda)\in \mathbb N$.}
\end{Remark}

\begin{Theorem}\label{c11} Let $x: [-1, 1]\to \mathbb R$ be the function whose graph is the intersection of the boundary $\partial K$ of $K$ with the half plane $\{(x,y)\in \mathbb R^2\; |\; x\geq 0\}$. Then $x(\cdot)$ has continuous
derivatives. In cases $A$ and $B$, the derivative $x^\prime$ is Lipschitz on $[-1,1]$, but in Case C it is not locally Lipschitz around $0$.
\end{Theorem}
\noindent{\bf Proof:} We first prove that $x^\prime$ exists and is continuous at $y=0$.
We use standard $(x,y)$ coordinates (for real and imaginary parts). Observe that $x(0)=1$. The concavity of $x(\cdot)$ implies
that for $y>0$ the slopes
\begin{equation*}
s(y):=\dfrac{x(y)-x(0)}{y}
\end{equation*}
have the property: $s(y_2)\geq s(y_1)$ if $y_2\leq y_1$. To calculate the limit of $s(y)$
as $y\to 0^+$, it is sufficient to choose a sequence $y_n\searrow 0$ and consider the limit
\begin{equation*}
s(0+):= \lim_{n\rightarrow\infty}\dfrac{x(y_n)-x(0)}{y_n}.
\end{equation*}
The same calculation for negative $y$ will result in the limit $s(0-)$. To conclude that $x$
is differentiable at $0$, we show that $s(0+)$ and $s(0-)$ both exist and equal 0. Note that $s(0-)=-s(0+)$.

Here is the calculation that establishes that $s(0+)=0$. Recall that
\begin{equation*}
T_n=\dfrac{e^{i\alpha_n}+e^{i\alpha_{n+1}}}{2}=\cos\left(\dfrac{\alpha_n-\alpha_{n+1}}{2}\right)e^{\displaystyle i\,\dfrac{\alpha_n+\alpha_{n+1}}{2}}
\end{equation*}
We now set
\begin{equation*}
y_n:= \mbox{\rm Im}\,(T_n) \quad \textrm{and} \quad x(y_n) := \mbox{\rm Re}\,(T_n)
\end{equation*}
and evaluate
\begin{equation*}
\lim_{n\to\infty}\dfrac{x(y_n)-x(0)}{y_n}= \lim_{n\to\infty}\dfrac{\cos\left(\dfrac{\alpha_n-\alpha_{n+1}}{2}\right) \cos\left(\dfrac{\alpha_n+\alpha_{n+1}}{2}\right) - 1}
{\cos\left(\dfrac{\alpha_n-\alpha_{n+1}}{2}\right) \sin\left(\dfrac{\alpha_n+\alpha_{n+1}}{2}\right)}
=0.
\end{equation*}
Thus, $x(\cdot)$ is differentiable at $y=0$ and $x^\prime(0)=0$. It follows that $x(\cdot)$ is differentiable on $[-1, 1]$, and $x^\prime$ is continuous away from the point $y=0$.

By the monotonicity of $x^\prime$ on $[-1, 1]$, the continuity of the derivative can be established by a similar argument. It is sufficient
to show that $x'(y_n)$ tends to zero as $n$ tends to infinity. We have
\begin{equation*}
x'(y_n) =
\dfrac{-\sin\left(\dfrac{\alpha_n+\alpha_{n+1}}{2}\right)}{\cos\left(\dfrac{\alpha_n+\alpha_{n+1}}{2}\right)}.
\label{eq:xprime-yn}
\end{equation*}
Again the limit is zero which proves the continuity of the derivative.

From Lemma \ref{lemma12}, we see that in cases $A$ and $B$ the sequences $\{r_n\}$ are bounded.
The curve $\partial K$ is given by
a linear function in the flat pieces which gives $x^{\prime\prime}(y)=0$, or by $x=x(y)$ where the second derivative of $x(\cdot)$ (except at the joints of the construction) is related to the curvature $1/r_n$ by
\begin{equation*}
\dfrac{1}{r_n}= \dfrac{x''}{(1+x'^2)^{3/2}}.
\end{equation*}

We need to prove that in cases $A$ and $B$, $x'(\cdot)$ is Lipschitz on $[-1,1]$.
As noted above, in these cases $x''(\cdot)$ exists (except at the joints) and is uniformly bounded on $[-1,1]$ by the facts that $\{r_n\}$ is bounded and $x'(\cdot)$ is continuous on $[-1,1]$. Thus, it is well-known that $x'$ is absolutely continuous on $[-1,1]$; see, e.g., \cite[Exercise 3.23, p.p.82]{G}. By Lebesgue's
Theorem (\cite[Theorem 6, Section 33]{KF}), we have
\begin{equation*}
 x'(y_2) -x'(y_1)= \int_{y_1}^{y_2} x''(s)\,ds,
\end{equation*}
where $y_1, y_2\in [-1,1]$.  By the bounded property of $x''(\cdot)$, the function $x'(\cdot)$ is Lipschitz on $[-1,1]$.

Note however that in Case $C$, the sequence $\{r_n\}$ tends to zero and therefore
$x''$ is unbounded in any neighborhood of $y=0$. This implies that in this case $x'$ is not locally Lipschitz
around $y=0$.
\qed

\begin{Remark} {\rm Since $x'(\cdot)$ is decreasing, $x(\cdot)$ is a concave function on $[-1,1]$. Equivalently, $-x(\cdot)$ is a convex on $[-1,1]$, and hence $x(\cdot)$ is locally Lipschitz on $(-1,1)$. Thus, we can apply \cite[Corollary 2.2.4, p.p.33]{CL} to obtain the continuity of $x'(\cdot)$ on $(-1,1)$, and hence on $[-1,1]$, from its differentiability on this interval. However, we give a direct proof as above for the convenience of the reader.
}\end{Remark}

\begin{Definition} For the set $K$ with the properties specified in Theorem \ref{c11}, we say that the $\partial K$ is $C^{1,1}$ around $(1,0)$ in cases $A$ and $B$, while $K$ has smooth boundary but $\partial K$ is not $C^{1,1}$ around $(1,0)$ in Case C.
\end{Definition}

\section{The Metric Projection}
\setcounter{equation}{0}

Given a nonempty closed convex set $\Omega\subset \mathbb{R}^m$, the metric projection from a given point
$x_0\in \mathbb{R}^m$ to $\Omega$ is defined by
\begin{equation*}
\Pi(x_0; \Omega):=\{w\in \Omega\; |\; d(x_0; \Omega)=\|x_0-w\|\},
\end{equation*}
where $d(x_0; \Omega):=\inf\{\|x_0-w\|\; |\; w\in \Omega\}$. It is well-known that $\Pi(x_0;\Omega)\in\Omega$ is always a singleton. Moreover, the mapping $\Pi(\cdot;\Omega)$ is nonexpansive in the sense that
\begin{equation*}
\|\Pi(x; \Omega)-\Pi(y; \Omega)\|\leq \|x-y\|\; \mbox{\rm for all }x,y\in \Bbb R^m.
\end{equation*}
The readers are referred to \cite{Urruty2001,bmn,r} for more details on the metric projection mapping.

In what follows, we consider the metric projection mapping $\Pi(\cdot; K)$, where the set $K$ is defined in the previous section. We omit $K$ in $\Pi(\cdot; K)$  if no confusion occurs.

The \emph{directional derivative} of the metric projection mapping at $x_0\not\in\Omega$ in the direction $v$ is given by
\begin{equation*}
D_v\Pi(x_0):=\lim_{t\to 0+}\,\dfrac{\Pi(x_0+tv)-\Pi(x_0)}{t}.
\label{eq:dir-deriv}
\end{equation*}

\begin{figure}[pbth]
\center
\includegraphics[width=0.4\linewidth]{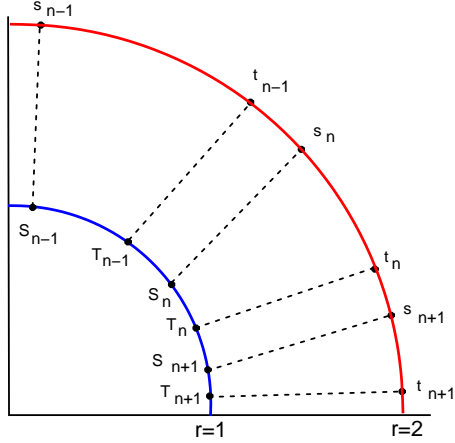}
\caption{ \emph{The construction of the projection of the convex set $K$.}}
\label{fig:convexset2}
\end{figure}

Now consider the parametrization of the circle $\mathcal{C}$ centered at the origin with radius 2: $x(\theta)=2e^{i\theta/2}$.

\begin{Lemma}\label{lm5} The directional derivative of $\Pi$ at $x(0)$ in the direction $v:=x^\prime(0)$ exists if and only if the limit
$$\lim_{\theta\to 0^+}\dfrac{\Pi(x(\theta))-\Pi(x(0))}{\theta-0}$$
exists.
\end{Lemma}
\noindent {\bf Proof:} By the nonexpansive property of the metric projection mapping, the following holds for any $\theta>0$:
\begin{align*}
\big \|\dfrac{\Pi(x(\theta))-\Pi(x(0))}{\theta-0}-\dfrac{\Pi(x(0)+\theta v)-\Pi(x(0))}{\theta-0}\big\|&=\big \|\dfrac{\Pi(x(\theta))-\Pi(x(0)+\theta v)}{\theta-0}\big\|\\
&\leq \big \|\dfrac{x(\theta)-x(0)-\theta v}{\theta-0}\big \|\\
&=\big \|\dfrac{x(\theta)-x(0)}{\theta-0}-v\big \|.
\end{align*}
Since $\lim_{\theta\to 0^+}\big \|\dfrac{x(\theta)-x(0)}{\theta-0}-v\big \|=0$, the conclusion follows easily. \qed

By Lemma \ref{lm5}, the directional derivative of the metric projection mapping at $(2,0)$ in the direction of the unit vector $i$ exists if and only if  $\dfrac{d}{d\theta}\Pi(x(\theta))|_{\theta=0}$ exists.

To better understand the metric projection mapping from the circle $\mathcal{C}$ onto $K$, we define
two points $2e^{it_n/2}$ and $2e^{is_n/2}$ such that
$$
\Pi(2e^{it_n/2})=T_n \quad \textrm{and} \quad \Pi(2e^{is_n/2})=S_n,
$$
where $T_n$ and $S_n$ are defined as before. The situation is depicted in Figure \ref{fig:convexset2}.


\begin{figure}[!ht]
\center
\includegraphics[width=0.6\linewidth]{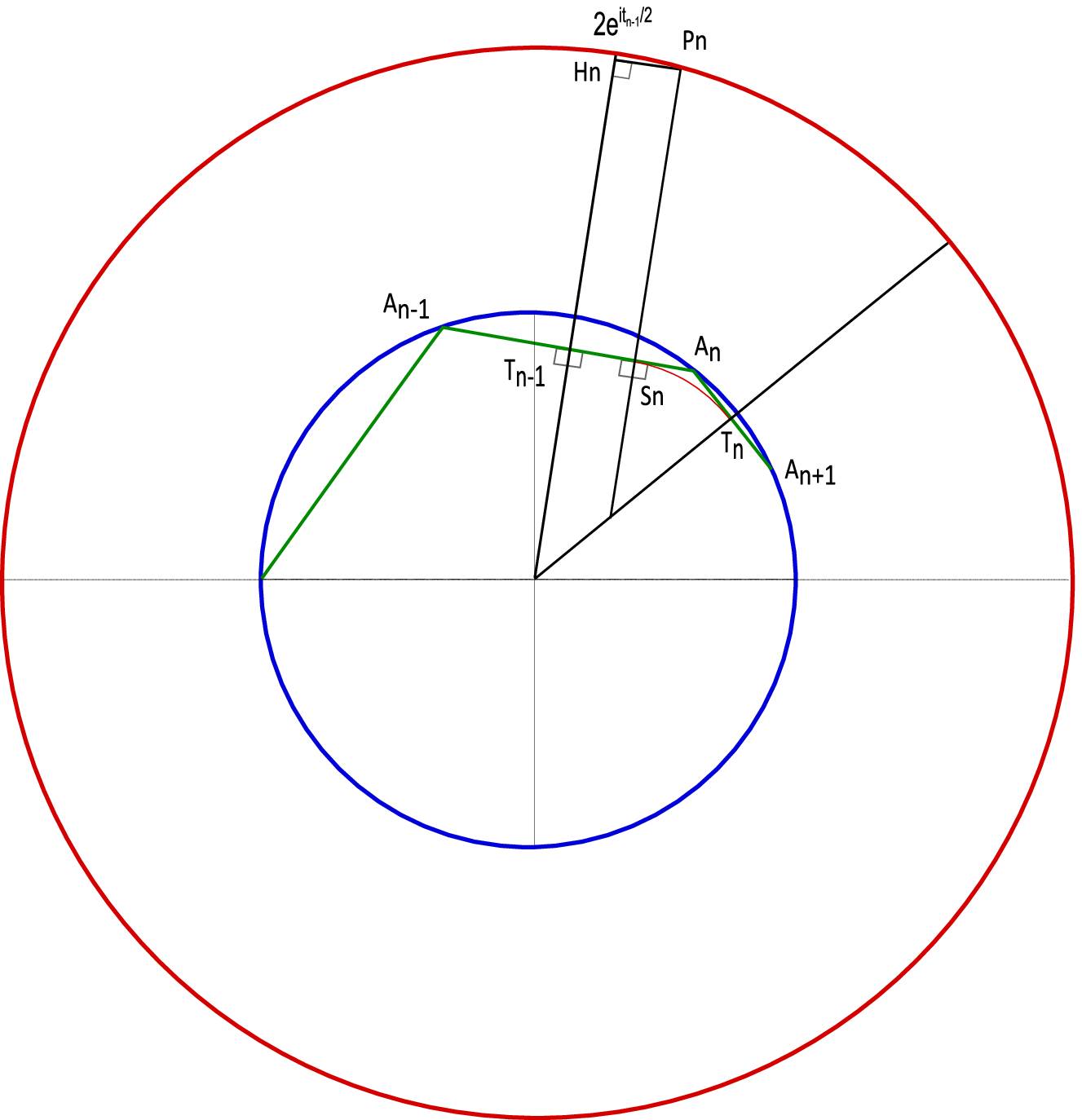}
\caption{An illustration for the proof of Lemma \ref{lm6}.}
\label{fig:proofoflm6}
\end{figure}

\begin{Lemma}\label{lm6} For any sequence $\{\alpha_n\}$ that defines our convex set $K$, we have
\begin{equation}\label{limit1}
\lim_{n\to\infty}\dfrac{\Pi(2e^{it_{n-1}/2})-\Pi(2e^{is_n/2})}{t_{n-1}-s_n}=i.
\end{equation}
\end{Lemma}
\noindent {\bf Proof:} Let
$$z_n:=\dfrac{\Pi(2e^{it_{n-1}/2})-\Pi(2e^{is_n/2})}{t_{n-1}-s_n}.$$
It suffices to show that
$$\|z_n\|\to 1\; \mbox{\rm and } \mbox{\rm arg}\,z_n\to \pi/2\; \mbox{\rm as }n\to \infty.$$
For the magnitude,  let $P_n$ denote point $2e^{is_n/2}$ and consider the orthogonal projection $H_n$ of $P_n$ onto the radii connecting the origin and $2e^{it_{n-1}/2}$ as seen in Figure \ref{fig:proofoflm6}. Obviously, $P_nH_n\Pi(2e^{it_{n-1}/2})\Pi(2e^{is_n/2})$ forms a rectangle. Opposite side lengths are equal, so $$\|\Pi(2e^{it_{n-1}/2})-\Pi(2e^{is_n/2})\| = \|P_n-H_n\|.$$
Considering the radii connecting the origin to the points $2e^{it_{n-1}/2}$ and $2e^{is_n/2}$ which mark off the angle $(t_{n-1}-s_n)/2$, we see that
\begin{equation*}
\|\Pi(2e^{it_{n-1}/2})-\Pi(2e^{is_n/2})\| = \|P_n-H_n\|=2\sin\left(\dfrac{t_{n-1}-s_n}{2}\right).
\end{equation*}
By the fundamental sine identity,
\begin{equation*}
\lim_{n\to\infty}\left\|\dfrac{\Pi(2e^{it_{n-1}/2})-\Pi(2e^{is_n/2})}{t_{n-1}-s_n}\right\| = 1.
\end{equation*}

To show that the argument tends to $\frac{\pi}{2}$, observe that since $T_{n-1}$ is the midpoint of $A_{n-1}A_n$, the line segment $A_{n-1}A_n$ is perpendicular to the line through $T_{n-1}$ and the origin. Since $\Pi(2e^{it_{n-1}/2})$ and $\Pi(2e^{is_n/2})$ are on the line segment $A_{n-1}A_n$ by definition, we get that
\begin{equation*}
\arg\left(\dfrac{\Pi(2e^{it_{n-1}/2})-\Pi(2e^{is_n/2})}{t_{n-1}-s_n}\right) = \frac{\pi}{2}+\arg\left(2e^{it_{n-1}/2}\right)=
\frac{\pi+t_{n-1}}{2}.
\end{equation*}
Observe that $2e^{it_{n-1}/2}, T_{n-1}$, and the origin are collinear (as in Figure \ref{fig:convexset1}), we have $t_{n-1}=\alpha_{n-1}+\alpha_n$. Thus,
\begin{equation*}
\lim_{n\to\infty}\arg\left(\dfrac{\Pi(2e^{it_{n-1}/2})-\Pi(2e^{is_n/2})}{t_{n-1}-s_n}\right) = \lim_{n\to\infty}\frac{\pi+\alpha_{n-1}+\alpha_n}{2} = \frac{\pi}{2}.
\end{equation*}
We have shown that the limit in \eqref{limit1} is $i$ as desired. \qed

Throughout the next few lemmas, we use $f(n) \sim g(n)$ to denote $\lim_{n\to\infty}f(n)/g(n)=1$.

 \begin{Lemma}\label{lme} If positive functions $f(n),g(n),h(n)$ satisfy $g(n)\sim h(n)$ and there exists a constant $c>0$ such that $\left |\dfrac{f(n)}{h(n)}-1\right|\geq c$ for all sufficiently large $n$, then $f(n)-g(n)\sim f(n)-h(n)$.
 \end{Lemma}
\noindent\textbf{Proof:} For all sufficiently large $n$, one has
\begin{equation*}
\left |\frac{f(n)-g(n)}{f(n)-h(n)}-1\right |=\left|\frac{g(n)-h(n)}{f(n)-h(n)}\right|=\left |\frac{\frac{g(n)}{h(n)}-1}{\frac{f(n)}{h(n)}-1}\right|\leq \frac{1}{c} \left |\frac{g(n)}{h(n)}-1\right|.
\end{equation*}
Then the conclusion follows easily. \qed

\begin{Lemma}\label{lme1} For any sequence $\{\alpha_n\}$ satisfying condition $(1.1)$, define
$$f(n):=\alpha_{n-1}-\alpha_{n+1},\; h(n):=\frac{\alpha_{n-1}-2\alpha_n+\alpha_{n+1}}{2}.$$
Then $f(n)$ and $h(n)$ satisfy the condition in Lemma \ref{lme}, i.e., exists a constant $c>0$ such that $\left|\dfrac{f(n)}{h(n)}-1\right|\geq c$ for all sufficiently large $n$.
\end{Lemma}
\noindent \textbf{Proof:} Define $b_n=\alpha_{n-1}-\alpha_n$. By condition $(1.1)$, $\{b_n\}$ is a positive decreasing sequencing that tends to $0$. Then $f(n)=b_n+b_{n+1}$ and $h(n)=\frac{b_{n}-b_{n+1}}{2}$. It suffices to show that there exists a constant $c>0$ such that $\left |\dfrac{2(b_n+b_{n+1})}{b_{n}-b_{n+1}}-1\right |\geq c$ for all sufficiently large $n$. Indeed,
$$\left |\dfrac{2(b_n+b_{n+1})}{b_{n}-b_{n+1}}-1\right |=\left |\dfrac{b_n+3b_{n+1}}{b_{n}-b_{n+1}}\right |\geq \dfrac{b_n+3b_{n+1}}{b_{n}}\geq 1\; \mbox{\rm for all }n\in \mathbb N.$$
The proof is now complete. \qed

\begin{Lemma}\label{lm7} For any sequence $\{\alpha_n\}$ that defines the convex set $K$, we have
$$\lim_{n\to\infty}\left(\dfrac{\Pi(2e^{is_n/2})-\Pi(2e^{it_n/2})}{s_n-t_n}-\dfrac{2(\alpha_n-\alpha_{n+1})}{\alpha_{n-1}+2\alpha_n-3\alpha_{n+1}}\cdot i\right)= 0.$$
\end{Lemma}
\noindent {\bf Proof:} Following the proof of Lemma \ref{lm6}, we compute the argument and magnitude separately.

Observe from the proof of Lemma \ref{lm1} that $A_nS_nO_nT_n$ is a right kite, and thus has perpendicular diagonals. In particular, this implies \begin{equation*}
\arg\left(\dfrac{\Pi(2e^{is_n/2})-\Pi(2e^{it_n/2})}{s_n-t_n}\right) = \frac{\pi}{2}+\arg\left(A_n-O_n\right) = \frac{\pi}{2}+\arg(T_n)+\dfrac{\pi-\psi_n}{2},
 \end{equation*}
where $\psi_n$ refers to the double-marked angle in Figure \ref{fig:convexset1}.

As noted from the proofs of Lemma \ref{lm1} and Theorem \ref{c11},
\begin{equation*}
\pi-\psi_n = \dfrac{\alpha_{n-1}-\alpha_{n+1}}{2}\; \mbox{\rm and }\arg\,(T_n) = t_n/2 = \dfrac{\alpha_n+\alpha_{n+1}}{2}.
\end{equation*}
Then
\begin{equation*}
\lim_{n\to\infty}\arg\left(\dfrac{\Pi(2e^{is_n/2})-\Pi(2e^{it_n/2})}{s_n-t_n}\right)=\lim_{n\to\infty}\left(\frac{\pi}{2}+
\frac{\alpha_n+\alpha_{n+1}}{2}+\dfrac{\alpha_{n-1}-\alpha_{n+1}}{4}\right)=\frac{\pi}{2}.
\end{equation*}
Now we compute the magnitude of the expression in question.  By formula \eqref{formularn}, as $S_n$ and $T_n$ are on the circle of radius $r_n$ centered at $O_n$, we see that
\begin{equation*}
\|\Pi(2e^{is_n/2})-\Pi(2e^{it_n/2})\| = 2r_n\sin\left(\dfrac{\pi-\psi_n}{2}\right) \sim r_n\cdot\dfrac{\alpha_{n-1}-\alpha_{n+1}}{2} \sim \alpha_n-\alpha_{n+1}.
\end{equation*}
We also have
\begin{equation*}
s_n-t_n=(t_{n-1}-t_n)-(t_{n-1}-s_n)= (\alpha_{n-1}-\alpha_{n+1}) -(t_{n-1}-s_n).
\end{equation*}

By Lemma \ref{lm6},
\begin{equation*}
t_{n-1}-s_n\sim \|\Pi(2e^{it_{n-1}/2})-\Pi(2e^{is_n/2})\|.
\end{equation*}

By the definition of $T_{n-1}=\Pi(2e^{it_{n-1}/2})$ and $S_n=\Pi(2e^{is_n/2})$, we see that
\begin{equation*}
\|\Pi(2e^{it_{n-1}/2})-\Pi(2e^{is_n/2})\| = \dfrac{\|A_n-A_{n-1}\|-\|A_{n+1}-A_n\|}{2},
\end{equation*}
 so that
\begin{equation*}
\|\Pi(2e^{it_{n-1}/2})-\Pi(2e^{is_n/2})\| = \sin\left(\dfrac{\alpha_{n-1}-\alpha_n}{2}\right)-\sin\left(\dfrac{\alpha_n-\alpha_{n+1}}{2}\right) \sim \dfrac{\alpha_{n-1}-2\alpha_n+\alpha_{n+1}}{2}.
\end{equation*}

Applying Lemma \ref{lme} and Lemma \ref{lme1} with
\begin{equation*}
f(n)=\alpha_{n-1}-\alpha_{n+1}, g(n) = t_{n-1}-s_n, h(n) = \dfrac{\alpha_{n-1}-2\alpha_n+\alpha_{n+1}}{2}
\end{equation*}
yields
\begin{equation*}
s_n-t_n \sim (\alpha_{n-1}-\alpha_{n+1}) - \dfrac{\alpha_{n-1}-2\alpha_n+\alpha_{n+1}}{2}=\dfrac{\alpha_{n-1}+2\alpha_n-3\alpha_{n+1}}{2}.
\end{equation*}

Then using the above three equations together, we get that
\begin{align*}
\lim_{n\to\infty}\left\|\dfrac{\Pi(2e^{is_n/2})-\Pi(2e^{it_n/2})}{s_n-t_n}\right\|=\lim_{n\to\infty}\dfrac{2(\alpha_n-\alpha_{n+1})}{\alpha_{n-1}+2\alpha_n-3\alpha_{n+1}}
\end{align*}
as desired. \qed


It is well-known the differentiability and the directional differentiability of the metric projection mapping are related to the second-order behavior of the boundary of the set involved; see \cite{Abat,FP,bor,Shap94} and the references therein. Note that the differentiability implies the directional differentiability. In the theorem below, we provide an example of a set with $C^{1,1}$ boundary but the metric projection mapping fails to be directionally differentiable.

\begin{Theorem}
In Case $B$, $\partial K$ is $C^{1,1}$ around $(1,0)$ and  $D_v\Pi$ does not exist at $x(0)=(2,0)$, where $v=x^\prime(0)=(0,1)$.

In Case $C$, $\partial K$ is $C^{1}$ but not $C^{1,1}$ around $(1,0)$, and $D_v\Pi$ does not exist at $x(0)=(2,0)$, where $v=x^\prime(0)=(0,1)$.
\end{Theorem}
\noindent {\bf Proof.} By Lemma \ref{lm5}, it suffices  to study the limit:
\begin{equation}
\lim_{\theta\rightarrow 0+} \dfrac{\Pi(2e^{i\theta/2})-\Pi(2e^{0})}{\theta}
\label{eq:the-limit}
\end{equation}
Let us first focus on Case B. Applying Lemma \ref{lm7}, we see that
\begin{equation}
\lim_{n\rightarrow \infty}\,\dfrac{\Pi(2e^{is_n/2})-\Pi(2e^{it_n/2})}{s_n-t_n}= \dfrac{2\lambda i}{3\lambda + 1}.
\label{eq:speeds1}
\end{equation}

By definition $T_n=\ell_ne^{\frac i2 (\alpha_n+\alpha_{n+1})}$ where $\ell_n=\cos\left(\dfrac{\alpha_n-\alpha_{n+1}}{2}\right)$ tends to $1$. Note that in Figure \ref{fig:convexset1} $T_n$, $O_n$, and the
origin are collinear. It follows that $t_n=\alpha_n+\alpha_{n+1}$. Since $2e^{it_n/2}$ projects to $T_n$,
we must have
\begin{equation*}
\lim_{n\rightarrow \infty}\dfrac{\Pi(2e^{it_n/2})-\Pi(2e^{0})}{t_n}=\dfrac i2
\label{eq:speeds2}
\end{equation*}
We write $\dfrac{\Pi(2e^{it_{n-1}/2})-\Pi(2e^{0})}{t_{n-1}}$ as a weighted mean of three fractions:
\begin{equation}
\dfrac{\Pi(2e^{it_{n-1}/2})-\Pi(2e^{is_n/2})}{t_{n-1}-s_n}\cdot \dfrac{t_{n-1}-s_n}{t_{n-1}}+
\dfrac{\Pi(2e^{is_n/2})-\Pi(2e^{it_n/2})}{s_n-t_n}\cdot \dfrac{s_n-t_n}{t_{n-1}}+
\dfrac{\Pi(2e^{it_n/2})-\Pi(2e^{0})}{t_n}\cdot \dfrac{t_n}{t_{n-1}}.
\label{eq:weighted-mean}
\end{equation}

Similarly, we write
\begin{equation}
\dfrac{\Pi(2e^{is_n/2})-\Pi(2e^{0})}{s_n}=
\dfrac{\Pi(2e^{is_n/2})-\Pi(2e^{it_n/2})}{s_n-t_n}\cdot \dfrac{s_n-t_n}{s_n}+
\dfrac{\Pi(2e^{it_n/2})-\Pi(2e^{0})}{t_n}\cdot \dfrac{t_n}{s_n}.
\label{eq:wm2}
\end{equation}

 Now we will show that the limit in  \eqref{eq:the-limit}, and hence the directional derivative of the metric projection mapping at $x(0)=(2,0)$ in the direction $v=(0,1)$, does not exist in case B. Suppose to the contrary that that this limit does exist.
Then
\begin{equation*}
\lim_{n\to\infty}\dfrac{\Pi(2e^{is_n/2})-\Pi(2e^{0})}{s_n}=
\lim_{n\to\infty}\dfrac{\Pi(2e^{it_n/2})-\Pi(2e^{0})}{t_n}=i/2.
\end{equation*}
Let $\lambda_n=\dfrac{s_n-t_n}{s_n}$ and $\beta_n=\dfrac{t_n}{s_n}$. Obviously, $\{\lambda_n\}$ and $\{\beta_n\}$ are nonnegative bounded sequences with
$$\lambda_n+\beta_n=1\; \mbox{\rm for all }n\in \mathbb N.$$

We will show that $\{\lambda_n\}$ converges to $0$. By a contradiction, suppose that this is not the case. Then there exist $\epsilon_0>0$ and a subsequence of $\{\lambda_{n_k}\}$ of $\{\lambda_n\}$ such that $\lambda_{n_k}\geq \epsilon_0$ for all $k\in \mathbb N$. By extracting a further convergent subsequence, we can assume without loss of generality that $\lim_{k\to \infty}\lambda_{n_k}=c>0$. From \eqref{eq:speeds1} and \eqref{eq:wm2}, one has
\begin{equation*}
\frac{i}{2}=c\frac{2\lambda i}{3\lambda + 1}+(1-c)\frac{i}{2},
\end{equation*}
which implies $$\dfrac{1}{2}=c\dfrac{2\lambda}{3\lambda + 1}+(1-c)\dfrac{1}{2}.$$
Since $\frac{2\lambda}{3\lambda + 1}<1/2$, one has $$\frac{1}{2}=c\frac{2\lambda}{3\lambda + 1}+(1-c)\frac{1}{2}<c/2+(1-c)/2=1/2,$$ a contradiction. We have shown that $\lim_{n\to\infty}\dfrac{s_n-t_n}{s_n}=0$, and hence  $\lim_{n\to\infty}\dfrac{t_n}{s_n}=1.$

Now, taking the limit as $n$ approaches infinity in  \eqref{eq:weighted-mean}, we get that
\begin{equation}
\dfrac{i}{2} = \lim_{n\to\infty}\left(i\cdot\dfrac{t_{n-1}-s_n}{t_{n-1}}+\dfrac{2\lambda i}{3\lambda+1}\cdot\dfrac{s_n-t_n}{t_{n-1}}+\dfrac{i}{2}\cdot\dfrac{t_n}{t_{n-1}}\right).
\label{eq:contraB}
\end{equation}
Of course, from $t_n=\alpha_n+\alpha_{n+1}$, in case B we must have \[\lim_{n\to\infty}\dfrac{t_n}{t_{n-1}}=\lambda.\]

Since $\lim_{n\to\infty}\dfrac{t_n}{s_n}=1$, we get \[\lim_{n\to\infty}\dfrac{s_{n}}{t_{n-1}}=\lambda.\]

Plugging these limits into  \eqref{eq:contraB} yields
\begin{equation*}
\dfrac{i}{2}=i(1-\lambda)+\dfrac{i}{2}\lambda,
\end{equation*}
which is absurd. Therefore, the limit from \eqref{eq:the-limit} does not exist, and hence in case B, $D_v\Pi$ does not exist at $x(0)=(2,0)$ in the direction $v=(0,1)$.

The proof showing that the limit does not exist in case C is analogous. Once more, suppose to the contrary that the limit from  \eqref{eq:the-limit} exists. We first claim that the limit exists only if
$$\lim_{n\to\infty}\frac{t_n}{s_n}=1.$$

Indeed, applying Lemma \ref{lm7} in Case C yields
\begin{equation*}
\lim_{n\to\infty}\dfrac{\Pi(2e^{is_n/2})-\Pi(2e^{it_n/2})}{s_n-t_n}=0.
\end{equation*}
Using \eqref{eq:wm2} and taking into account that
\begin{equation*}
\lim_{n\to\infty} \left\|\dfrac{\Pi(2e^{is_n/2})-\Pi(2e^{0})}{s_n}\right\| = \dfrac{1}{2}\;\mbox{\rm and }\lim_{n\to\infty}\left\|\dfrac{\Pi(2e^{it_n/2})-\Pi(2e^{0})}{t_n}\right\|=\dfrac{1}{2},
\end{equation*}
one has
$$\lim_{n\to\infty}\dfrac{t_n}{s_n}=1.$$

From $t_n=\alpha_n+\alpha_{n+1}$, in Case C we must have \[\lim_{n\to\infty}\dfrac{t_n}{t_{n-1}}=0.\]

Since $\lim_{n\to\infty}\dfrac{t_n}{s_n}=1$, we get \[\lim_{n\to\infty}\dfrac{s_{n}}{t_{n-1}}=0.\]

Plugging these limits into  \eqref{eq:contraB} yields

\[\dfrac{i}{2}=i,\]
which is contradiction. Thus, the limit from  \eqref{eq:the-limit} does not exist in Case C as well, and hence $D_v\Pi$ does not exist at $x(0)=(2,0)$ in the direction $v=(0,1)$.   \qed

\begin{Remark} We conjecture that in Case $A$, $D_v\Pi$ does exist at $x(0)=(2,0)$ in the direction $v=(0,1)$.
\end{Remark}

{\bf Acknowledgements.}  The research of Nguyen Mau Nam was
partially supported by the NSF under grant \#1411817 and the Simons Foundation under grant \#208785. The research of J.J.P. Veerman was partially supported by the European Union's
Seventh Framework Program (FP7-BEGPOT-2012-2013-1) under grant agreement n316165.

\end{document}